\documentclass[11pt]{article}
\usepackage{dsfont} 
\usepackage{amssymb} 

\usepackage{mathrsfs} 

\begin{document}
\title{{\bf\Large Geometric shape of invariant manifolds\\ for a class of stochastic partial
differential equations
\thanks{This work was supported by the
National Science Foundation of China (Grants No. 10901115 and 11071177) and the NSF Grant  1025422.  } } }
\author{Guanggan Chen\\
College of Mathematics and Software Science,\\
Sichuan Normal University, Chengdu, 610068, China\\
E-mail: chenguanggan@hotmail.com\\
\\
Jinqiao Duan\\
Department of Applied Mathematics,\\
Illinois Institute of Technology, Chicago, IL 60616, USA\\
E-mail: duan@iit.edu\\
\\
Jian Zhang\\
College of Mathematics and Software Science,\\
Sichuan Normal University, Chengdu, 610068, China\\
E-mail: zhangjiancdv@sina.com }

\date{\today}
\maketitle
\baselineskip 8mm
\begin{center}
\begin{minipage}{100mm}
{\bf \large Abstract:} {\small Invariant manifolds play an important
role in the study of the qualitative dynamical behaviors for
nonlinear stochastic partial differential equations. However, the
geometric shape of these manifolds is largely unclear. The purpose
of the present paper is to try to describe the geometric shape of
invariant manifolds for  a class of stochastic partial differential
equations with multiplicative white noises.  The local geometric
shape of invariant manifolds is approximated, which holds with
significant likelihood. Furthermore, the result is compared with
that for the corresponding deterministic partial differential
equations.}

\end{minipage}
\end{center}

\renewcommand{\theequation}{\thesection.\arabic{equation}}
\setcounter{equation}{0}


\vspace{1cm}

\section{ Introduction }

\quad\, Stochastic partial differential equations arise as models
for various complex systems under random influences. There have been
recent rapid     progresses in this area (see \cite{B, DZ, GS, R,
WD}). For stochastic partial differential equations, random
invariant manifolds play an important role in the study of dynamics
because they provide a geometric structure to understand or reduce
stochastic dynamics. Although the existence for such random
invariant manifolds is established for certain stochastic partial
differential equations (e.g.,  \cite{CDKS, CCL05, DLS03, DLS04,
LS07, MZZ08}), the geometric shape of these manifolds is largely
unclear. The purpose of the present paper is to try to describe the
geometric shape of invariant manifolds for a class of stochastic
partial differential equations.
\par
We consider a class of stochastic partial differential equations in
the following form
\begin{equation}
\frac{d u}{d t}+Lu-u^p=\sigma u\circ \dot{W},
\end{equation}
subject to the homogeneous Dirichlet boundary conditions on a
bounded domain with scalar   white noise $\dot{W}$  of
Stratonovich type. The linear operator $-L$ generates a
$C_0$-semigroup, which is given in detail in the next section. The
nonlinear power exponent $p$ belongs to $ (1, +\infty)$,  and
$\sigma$ is a real parameter in $(0, +\infty)$.
\par
It is well known that the theory of invariant manifolds has been
developed well for deterministic dynamical systems. However, for the
stochastic dynamical systems generated by stochastic partial
differential equations, due to their nonclassical fluctuation of
driving noise and infinite dimensionality, the theory of invariant
manifolds, together with their approximation and computation, is
still in its infancy.
\par

For stochastic partial differential equation, a random invariant
manifold has various samples in an infinite dimensional space.
Therefore it is difficult in general to describe or ``visualize"
random invariant manifolds, let alone the reduction of dynamics on
them. Blomker and Wang \cite{BW}, and Sun et al \cite{SDL} have done
some work on describing  such invariant manifolds. In this paper,
we will consider an approximate local geometric shape of invariant
manifolds for Equation (1.1).
\par
More precisely, for Equation (1.1), we first construct a local
invariant manifold. Then by approximating the invariant manifold
step by step, we establish an approximate local geometric shape of
the invariant manifold, which holds with probabilistic significance.
Next, we study the corresponding deterministic system of Equation
(1.1) (i.e., Equation (1.1) with $\sigma=0$). Using the same method,
we drive the invariant manifold and its approximating local
geometric shape, which always holds.

\par
This paper is organized as follows. In the  next section, we present
the assumptions of the linear operator $L$, introduce the basic
concepts on random dynamical systems and the random evolutionary
equation induced by Equation (1.1). In the third section, we show
Theorem 3.1 on the existence of the local random invariant manifold
for Equation (1.1). In the fourth section, we prove Theorem 4.1 on
the local geometric shape of the random invariant manifold.
Furthermore, we give an example to explain the local geometric shape
in Remark 4.1. In the fifth section, we discuss the local geometric
shape of the invariant manifold  for the corresponding deterministic
system of Equation (1.1). We comment on the results in the final
section. We consider only unstable invariant manifolds, as stable
invariant manifolds may be discussed similarly.
\par

\renewcommand{\theequation}{\thesection.\arabic{equation}}
\setcounter{equation}{0}

\section{Preliminaries}

{\bf 2.1.\quad Assumption of the linear operator $L$}
\par
Let $E$ be a separable Hilbert space with norm $\|\cdot\|$ and
scalar product $\langle \cdot, \cdot\rangle$, and $L$ be a closed
self-adjoint linear operator with dense domain $D(L)$ in $E$. Let
$id$ be the identity operator on $E$.
\par
{\bf Hypothesis}\quad{\it There exists a constant $a\geq 0$ such
that $(L+a \cdot\hbox{id})$ is positive and $(L+a
\cdot\hbox{id})^{-1}$ is compact.}
\par
This assumption implies that the spectrum of $L$ consists of only
eigenvalues with finite multiplicities,
\begin{equation}
-a<\lambda_1\leq \lambda_2\leq \cdots, \quad \lim\limits_{n\to
+\infty}\lambda_n=+\infty,
\end{equation}
and the associated eigenfunctions $\{e_n\}_{n\in \mathbb{N}}$,
$e_i\in D(L)\subset E$ form an othonormal basis of $E$. An example
of the linear operator $L$ is   $L=-\partial_{xx}-3\cdot id$ on
$H^1_0([0, \pi])$, whose eigenvalues are $\lambda_k=k^2-3$ with
the corresponding eigenfunctions $e_k=\sin kx$ $, k=1, 2, 3,\cdots$.
\par
Furthermore, the positivity of $(L+a \cdot\hbox{id})$ allows one to
define the fractional power of $(L+a \cdot\hbox{id})$, which we
denote by $(L+a \cdot\hbox{id})^\alpha$ for $\alpha\in [0, 1)$, see
Henry \cite{Henry} or  Temam \cite{T}. The domain of $(L+a
\cdot\hbox{id})^\alpha$, which we denote by $E^\alpha$, is a Hilbert
space with the scalar product $\langle u, \tilde{u}\rangle_\alpha=\langle
(L+a \cdot\hbox{id})^\alpha u, (L+a
\cdot\hbox{id})^\alpha\tilde{u}\rangle$ and corresponding norm
$|\cdot|_\alpha$.
\par

From (2.1), there exists $\lambda_N< 0$ such that $-a<\lambda_1\leq
\lambda_2\leq \cdots\leq \lambda_N< 0$ and $\lambda_{N+1}\geq 0$.
Denote $\lambda_u:=\lambda_N \;(< 0)$ and
$\lambda_s:=\lambda_{N+1}\; (\geq 0)$. Put $E_u:=span\{e_1,\cdots,
e_N\}$. Let $P_u$ be the orthogonal projection from $E$ to $E_u$ and
$P_s=I-P_u$. Put $L_u=P_u{L}$ and $L_s=P_s{L}$. In the following, we
use the subscript ``$u$" always for projection onto $E_u$ and the
subscript ``$s$" for projection onto $E_s$. Then $E=E_u\oplus E_s$
and $E^\alpha=E_u\oplus E_s^\alpha$, where $E_s^\alpha=E_s\bigcap
E^\alpha$ and $E_u\subset E^\alpha$ with $\alpha\in [0, 1)$.
\par
From Henry \cite{Henry}, there exists $M>0$ such that
\begin{equation}
\begin{array}{l}
\|e^{-L_st}P_s\|_{\mathcal{L}(E^\alpha, E^\alpha)}\leq
Me^{-\lambda_st},\quad
t\geq 0;\\
\|e^{-L_st}P_s\|_{\mathcal{L}(E^\alpha, E)}\leq
\frac{M}{t^\alpha}e^{-\lambda_st},\quad
t\geq 0;\\
\|e^{-L_ut}P_u\|_{\mathcal{L}(E^\alpha, E^\alpha)}\leq Me^{-\lambda_ut},\quad t\leq 0;\\
\|e^{-L_ut}P_u\|_{\mathcal{L}(E^\alpha, E)}\leq Me^{-\lambda_ut},
\quad t\leq 0,
\end{array}
\end{equation}
where $\mathcal{L}(X, Y)$ is the usual space of bounded linear
operator from Banach space $X$ to Banach space $Y$.
\par

\noindent {\bf 2.2.\quad Random dynamical systems}
\par

Let us recall some basic concepts in random dynamical systems as in \cite{DLS03}. Let
$(\Omega,\mathcal{F}, \mathds{P})$ be a probability space. A flow
$\theta$ of mappings $\{\theta_t\}_{t\in \mathbb{R}}$ is defined on
the sample space $\Omega$ such that
\begin{equation}
\theta: \mathbb{R}\times \Omega\to \Omega, \quad
\theta_0=id,\quad \theta_{t_1}\theta_{t_2}=\theta_{t_1+t_2},
\end{equation}
for $t_1, t_2 \in \mathbb{R}$. This flow is supposed to be
$(\mathcal{B}(\mathbb{R})\otimes\mathcal{F},
\mathcal{F})$-measurable, where $\mathcal{B}(\mathbb{R})$ is the
$\sigma$-algebra of Borel sets on the real line $\mathbb{R}$. To
have this measurability, it is not allowed to replace $\mathcal{F}$
by its $\mathds{P}$-completion $\mathcal{F}^{\mathds{P}}$; see
Arnold \cite{Arnold} p. 547. In addition, the measure $\mathds{P}$
is assumed to be ergodic with respect to $\{\theta_t\}_{t\in
\mathbb{R}}$. Then $(\Omega,\mathcal{F}, \mathds{P},\mathbb{R},
\theta)$ is called a metric dynamical system.
\par
For our applications, we will consider a special but very important
metric dynamical system induced by the Brownian motion. Let $W(t)$
be a two-sided Wiener process with trajectories in the space
$C_0(\mathbb{R},\mathbb{R})$ of real continuous functions defined on
$\mathbb{R}$, taking zero value at $t = 0$. This set is equipped
with the compact open topology. On this set we consider the
measurable flow $\theta = \{\theta_t\}_{t\in \mathbb{R}}$, defined
by $\theta_t\omega =\omega(\cdot+t)-\omega(t)$. The distribution of
this process generates a measure on
$\mathcal{B}(C_0(\mathbb{R},\mathbb{R}))$ which is called the Wiener
measure. Note that this measure is ergodic with respect to the above
flow; see the Appendix in Arnold [1]. Later on we will consider,
instead of the whole $C_0(\mathbb{R},\mathbb{R})$, a
$\{\theta_t\}_{t\in \mathbb{R}}$-invariant subset $\Omega\subset
C_0(\mathbb{R},\mathbb{R})$) of $\mathds{P}$-measure one and the
trace $\sigma$-algebra $\mathcal{F}$ of
$\mathcal{B}(C_0(\mathbb{R},\mathbb{R}))$ with respect to $\Omega$.
A set $\Omega$ is called {$\{\theta_t\}_{t\in \mathbb{R}}$-invariant
if $\theta_t\Omega = \Omega$ for $t \in \mathbb{R}$. On
$\mathcal{F}$, we consider the restriction of the Wiener measure
also denoted by $\mathds{P}$.
\par
The dynamics of the system on the state space $E$ over the flow
$\theta$ is described by a cocycle. For our applications it is
sufficient to assume that $(E, dE)$ is a complete metric space. A
cocycle $\phi$ is a mapping:
$$
\phi : \mathbb{R}^+\times\Omega\times E\to E,
$$
which is $(\mathcal{B}(\mathbb{R})\otimes \mathcal{F} \otimes
\mathcal{B}(E), \mathcal{F})$-measurable such that
$$
\begin{array}{l}
\phi(0, \omega, x) = x \in E,\\
\phi(t_1 + t_2, \omega, x) = \phi(t_2, \theta_{t_1}\omega, \phi(t_1,
\omega, x)),
\end{array}
$$
for $t_1, t_2 \in \mathbb{R}^+, \omega \in \Omega$ and $x\in E$.
Then $\phi$ together with the metric dynamical system $\theta$ forms
a random dynamical system.
\par

\par

\noindent {\bf 2.3.\quad Random evolutionary equation }
\par
We consider a linear stochastic differential equation
\begin{equation}
dz+zdt=\sigma dW.
\end{equation}
A solution of this equation is called an Ornstein-Uhlenbeck process.
We have the following results, see Duan, Lu and Schmalfuss \cite{DLS03, DLS04}.
\par
{\bf Lemma 2.1}\quad {\it (i)\quad There exists a
$\{\theta_t\}_{t\in \mathbb{R}}$-invariant set $\Omega\in
\mathcal{B}(C_0(\mathbb{R}, \mathbb{R}))$ of full measure with
sublinear growth
$$
\lim\limits_{t\to \pm\infty}\frac{|\omega(t)|}{|t|}=0,\quad
\omega\in \Omega
$$
of $\mathds{P}$-measure one.
\par
(ii)\quad For $\omega\in \Omega$, the random variable
$$
z(\omega)=-\sigma\int_{-\infty}^{0}e^\tau\omega(\tau)d\tau
$$
exists and generates a unique stationary solution of Equation (2.4)
given by
$$
z(\theta_t\omega)=-\sigma\int_{-\infty}^{0}e^\tau\theta_t\omega(\tau)d\tau
=-\sigma\int_{-\infty}^{0}e^\tau\omega(\tau+t)d\tau+
\sigma\omega(t).
$$
The mapping $t\to z(\theta_t\omega)$ is continuous.
\par
(iii)\quad In particular,
$$
\lim\limits_{t\to \pm\infty}\frac{|z(\theta_t\omega)|}{|t|}=0, \quad
for\quad \omega\in \Omega.
$$
\par
(iv)\quad In addition,
$$
\lim\limits_{t\to
\pm\infty}\frac{1}{t}\int_0^tz(\theta_\tau\omega)d\tau=0, \quad
for\quad \omega\in \Omega.
$$
}
\par
We now replace $\mathcal{B}(C_0(\mathbb{R},\mathbb{R}))$ by
$
\mathcal{F}=\{\Omega\bigcap F|\quad F\in
\mathcal{B}(C_0(\mathbb{R},\mathbb{R}))\}
$
for $\Omega$ given in Lemma 2.1. The probability measure is the
restriction of the Wiener measure to this new $\sigma$-algebra,
which is also denoted by $\mathds{P}$. In the following we will
consider the metric dynamical system
$
(\Omega,\mathcal{F}, \mathds{P}, \mathbb{R}, \theta).
$
\par
Now we show that the solution of Equation (1.1) defines a random
dynamical. Firstly, the equivalent It$\hat{o}$ equation of Equation (1.1) is given by
\begin{equation}
du=-Ludt+u^pdt+\frac{u}{2}dt+\sigma udW
\end{equation}
with the initial data $u(0)=u_0\in E^\alpha$ being $\mathcal{F}_0$-measurable.
Equation (2.5) can be written in the following mild integral form
$$
u(t)=e^{-Lt}u_0+\int_0^te^{-L(t-r)}(u^p(r)+\frac{u(r)}{2})dr+\sigma\int_0^te^{-L(t-r)}u(r)dW(r)
$$
almost surely for arbitrary $u_0\in E^\alpha$, in which the
stochastic integral is to interpret in the sense of It$\hat{o}$.
\par
Under the transformation of Ornstein-Uhlenbeck process (2.4),
Equation (2.5) becomes a random evolutionary equation (i.e., an
evolutionary equation with random coefficients)
\begin{equation}
\frac{dv}{dt}=-Lv +zv +e^{-z}F(e^zv)
\end{equation}
with $v(0)=u_0e^{-z(0)}:=x\in E^\alpha$, where $v=ue^{-z}$ and $F(v)=v^p$ with
$z=z(t):=z(\theta_t\omega)$. In contrast to the original stochastic
differential equation (1.1), no stochastic integral appears here. Then the
mild integral form of (2.6) is
$$
v(t)=e^{-Lt+\int_0^tz(\tau)d\tau}x+\int_0^te^{-L(t-r)+\int_r^tz(\tau)
d\tau}e^{-z(r)}F(e^{z(r)}v(r))dr
$$
almost surely for any $x\in E^\alpha$.
\par
Since our purpose is to consider the dynamical behavior of solution
of Equation (2.6) in a neighborhood of the fixed point $v=0$ in this
paper, now we introduce a truncated equation of Equation (2.6) such
that its nonlinear term has a small Lipschitz constant.
\par
Let $\chi: E^\alpha\to \mathbb{R}$
be a $C_0^\infty$ function, a cut-off function, such that
$$
\chi(v)= \left\{
\begin{array}{l}
1,\quad \hbox{if}\quad |v|_\alpha \leq 1, \\
0,\quad \hbox{if}\quad |v|_\alpha \geq 2.
\end{array}
\right.
$$
For any positive parameter $R$, we define
$\chi_R(v)=\chi(\frac{v}{R})$ for all $v\in E^\alpha$. Let
$F^{(R)}(v)=\chi_R(v)F(v)$.
For every $l_F>0$ and every $\omega\in
\Omega$, there must exist a positive random variable $R$ such that
\begin{equation}
\|F^{(R)}(v)-F^{(R)}(\tilde{v})\|\leq l_F| v-\tilde{v}|_\alpha.
\end{equation}
Then the truncated equation of Equation (2.6) is as
follows
\begin{equation}
\frac{dv}{dt}=-Lv +zv +e^{-z}F^{(R)}(e^zv).
\end{equation}
\par
By the classical evolutionary equation theory, Equation (2.8) has a
unique solution for every $\omega\in \Omega$. No exceptional sets
with respect to the initial conditions appear. Hence the solution
mapping
$$
(t, \omega,  x) \mapsto \phi(t, \omega)x:=v(t,\omega; x)
$$
generates a continuous random dynamical system. Indeed, the mapping
$\phi$ is $(\mathcal{B}(\mathbb{R})\otimes\mathcal{F}\otimes
\mathcal{B}(E^\alpha), \mathcal{F})$-measurable.
\par
Introduce the transform
$$
T(\omega, x)=xe^{-z(\omega)}
$$
and its inverse transform
$$
T^{-1}(\omega, x)=xe^{z(\omega)}
$$
for $x\in E^\alpha$ and $\omega\in \Omega$. Then for the random
dynamical system $v(t,\omega; x)$ generated by Equation (2.6),
$$
(t, \omega,  x) \mapsto T^{-1}(\theta_t\omega, v(t, \omega; T(w,
x))):=u(t, \omega; x)
$$
is the random dynamical system generated by Equation (1.1). For more
about the relation between (1.1) and (2.6), we refer to Duan, Lu and
Schmalfuss \cite{DLS03}.
\par

\par
\renewcommand{\theequation}{\thesection.\arabic{equation}}
\setcounter{equation}{0}

\section{Existence of local invariant manifolds}

\quad\quad In this section, we shall use the method of Duan, Lu and
Schmalfuss \cite{DLS04} to establish the local invariant manifold of
Equation (1.1).
\par
Define a Banach space for each $\beta\in (\lambda_u, \lambda_s)$ as
follows
$$
C_{\beta}^-=\{f(\cdot)\in C((-\infty, 0]; E^\alpha)|\quad
\sup\limits_{t\leq 0}e^{\beta t-\int_0^tz(\tau)d\tau}
|f|_\alpha<\infty \}
$$
with the norm
$$
\|f\|_{C_{\beta}^-}=\sup\limits_{t\leq 0}e^{\beta
t-\int_0^tz(\tau)d\tau} |f|_\alpha.
$$
\par
Since that
$$
\begin{array}{ll}
|e^{-L_s(t-r)+\int_r^tz(\tau)d\tau}P_su(r)|_\alpha &\leq
Me^{-\beta r+\int_0^rz(\tau)d\tau}e^{-\lambda_s(t-r)+\int_r^tz(\tau)d\tau}|u(r)|_{C_{\beta}^-}\\
&\leq
Me^{(\lambda_s-\beta) r}e^{-\lambda_s t+\int_0^tz(\tau)d\tau}|u(r)|_{C_{\beta}^-}\\
&\longrightarrow 0, \quad \hbox{as}\quad r\to -\infty,
\end{array}
$$
we have that
$$
P_sv(t)=\int_{-\infty}^t
e^{-L_s(t-r)+\int_r^tz(\tau)d\tau}e^{-z(r)}F_s^{(R)}(e^{z(r)}v(r))
dr.
$$
Then we have the following result. For the detailed proof, please see
Duan, Lu and Schmalfuss \cite{DLS04}.
\par
{\bf Lemma 3.1}\quad{\it Suppose that $v(\cdot)$ is in
$C_{\beta}^-$. Then $v(t)$ is the solution of Equation (2.8) with
the initial datum $v(0)=x$ if and only if $v(t)$ satisfies
\begin{equation}
\begin{array}{lll}
v(t)&=&e^{-L_ut+\int_0^tz(\tau)
d\tau}\xi+\int_0^te^{-L_u(t-r)+\int_r^tz(\tau)d\tau}e^{-z(r)}F_u^{(R)}(e^{z(r)}v(r))dr\\
&&+\int_{-\infty}^t
e^{-L_s(t-r)+\int_r^tz(\tau)d\tau}e^{-z(r)}F_s^{(R)}(e^{z(r)}v(r))
dr,
\end{array}
\end{equation}
where $\xi=P_ux\in E_u$. }
\par
Define
\begin{equation}
\begin{array}{lll}
J(v, \xi)&=&e^{-L_ut+\int_0^tz(\tau)
d\tau}\xi+\int_0^te^{-L_u(t-r)+\int_r^tz(\tau)d\tau}e^{-z(r)}F_u^{(R)}(e^{z(r)}v(r))dr\\
&&+\int_{-\infty}^t
e^{-L_s(t-r)+\int_r^tz(\tau)d\tau}e^{-z(r)}F_s^{(R)}(e^{z(r)}v(r))
dr,
\end{array}
\end{equation}
and also denote
\begin{equation}
SC:=Ml_F[\frac{1}{\beta-\lambda_u}+\frac{\Gamma(1-\alpha)}{(\lambda_s-\beta)^{1-\alpha}}], \label{SC}
\end{equation}
where $M$ is the positive constant in (2.2), $l_F$ is the Lipschitz
constant in (2.7), $\Gamma(\cdot)$ is the Gamma function, $\alpha\in
[0, 1)$ and $\beta\in (\lambda_u, \lambda_s)$.
\par
Then
\begin{equation}
\begin{array}{lll}
\|J(v, \xi)-J(\tilde{v}, \xi)\|_{C_\beta^-}&\leq&
\|\int_0^te^{-L_u(t-r)+\int_r^tz(\tau)d\tau}e^{-z(r)}[F_u^{(R)}(e^{z(r)}v(r))-F_u^{(R)}(e^{z(r)}\tilde{v}(r))]dr\\
&&+\int_{-\infty}^t
e^{-L_s(t-r)+\int_r^tz(\tau)d\tau}e^{-z(r)}[F_s^{(R)}(e^{z(r)}v(r))-F_s^{(R)}(e^{z(r)}\tilde{v}(r))]
dr\|_{C_\beta^-}\\
&\leq & Ml_F\cdot \sup\limits_{t\leq 0}[\int_0^te^{(\beta-\lambda_u)(t-r)}dr+
\int_{-\infty}^t\frac{1}{(t-r)^\alpha}e^{(\beta-\lambda_s)(t-r)}dr]\|v-\tilde{v}\|_{C_\beta^-}\\
&\leq & SC \|v-\tilde{v}\|_{C_\beta^-}.
\end{array}
\end{equation}
\par
Let $SC<1$. Then by the uniform contraction mapping principle, for each $\xi\in  E_u$, $J(v,\xi)$ has
a unique fixed point $v^*(t, \omega; \xi)\in C_{\beta}^-$. Put $h(\omega, \xi)=P_sv^*(0, \omega; \xi)$. Thus
\begin{equation}
\begin{array}{lll}
h(\omega, \xi)=\int_{-\infty}^0
e^{L_sr+\int_r^0z(\tau)d\tau}e^{-z(r)}F_s^{(R)}(e^{z(r)}v(r)) dr.
\end{array}
\end{equation}
\par
{\bf Lemma 3.2}\quad{\it Let $R$ be a positive random variable such
that $l_F$ satisfies $SC<1$. For the unique fixed point $v^*=v^*(t,
\omega; \xi)=J(v^*)\in C_\beta^-$ of the operator $J$, there exists a
positive constant $C$ such that
$$
\|v^*(t,\omega;\xi_1)-v^*(t,\omega;\xi_2)\|_{C_\beta^-}\leq
C|\xi_1-\xi_2|_\alpha.
$$
Moreover,
$$
\|h(\omega,\xi_1)-h(\omega,\xi_2)\|_{C_\beta^-}\leq
C|\xi_1-\xi_2|_\alpha.
$$
}
\par
{\bf Lemma 3.3} \quad{\it Let $R$ be a positive random variable such
that $l_F$ satisfies $SC<1$. Then there exists a positive constant
$C$ such that
\begin{equation}
\begin{array}{l}
\|v^*(t,\omega;\xi)\|_{C_\beta^-}\leq C|\xi|_\alpha,\\
\|v_s^*(t,\omega;\xi)\|_{C_\beta^-}\leq C|\xi|_\alpha,\\
\|v_u^*(t,\omega;\xi)\|_{C_\beta^-}\leq C|\xi|_\alpha,
\end{array}
\end{equation}
where $v_s^*=P_sv^*$ and $v_u^*=P_uv^*$. }
\par
{\bf Proof.}\quad Firstly, for $t\leq 0$, since $\beta\in
(\lambda_u, \lambda_s)$, we have
\begin{equation}
\begin{array}{ll}
\|J(0, \xi)\|_{C_\beta^-}&=\|e^{-L_ut+\int_0^tz(\tau)
d\tau}\xi\|_{C_\beta^-}\leq \sup\limits_{t\in (-\infty, 0]}e^{\beta
t-\int_0^tz(\tau)d\tau}e^{-L_ut+\int_0^tz(\tau) d\tau}|\xi|_\alpha\\
&\leq \sup\limits_{t\in (-\infty, 0]}M e^{(\beta-\lambda_u)
t}|\xi|_\alpha\leq C|\xi|_\alpha.
\end{array}
\end{equation}
It follows from Lemma 3.2, (3.4) and (3.7) that
$$
\|v^*(t,\omega;\xi)\|_{C_\beta^-}\leq \|J(v^*,\xi)-J(0,
\xi)\|_{C_\beta^-}+\|J(0, \xi)\|_{C_\beta^-}\leq
SC\|v^*(t,\xi)\|_{C_\beta^-}+C|\xi|_\alpha,
$$
which implies that $\|v^*(t,\omega;\xi)\|_{C_\beta^-}\leq
\frac{C}{1-SC}|\xi|_\alpha$.
\par
Meanwhile,
$$
\begin{array}{ll}
\|v_s^*(t,\omega;\xi)\|_{C_\beta^-}&=\|P_sv^*(t,\omega;\xi)\|_{C_\beta^-}\\
&=\|\int_{-\infty}^t
e^{-L_s(t-r)+\int_r^tz(\tau)d\tau}e^{-z(r)}F_s^{(R)}(e^{z(r)}v(r))
dr\|_{C_\beta^-}\\
&\leq \int_{-\infty}^t e^{\beta
t-\int_0^tz(\tau)d\tau}\frac{M}{(t-r)^\alpha}e^{-\lambda_s(t-r)}e^{\int_r^tz(\tau)d\tau}
e^{-\beta
r+\int_0^rz(\tau)d\tau}l_F\|v(r)\|_{C_\beta^-} dr\\
&\leq Ml_F\|v(r)\|_{C_\beta^-} \int_{-\infty}^t
e^{(\beta-\lambda_s)(t-r)
}\frac{1}{(t-r)^\alpha}dr\\
&\leq C\|v(r)\|_{C_\beta^-}\\
&\leq C|\xi|_\alpha.
\end{array}
$$
\par
Therefore
$$
\|v_u^*(t,\omega;\xi)\|_{C_\beta^-}=
\|v^*(t,\omega;\xi)-v_s^*(t,\omega;\xi)\|_{C_\beta^-}\leq
\|v^*(t,\omega;\xi)\|_{C_\beta^-}+\|v_s^*(t,\omega;\xi)\|_{C_\beta^-}\leq
C|\xi|_\alpha.
$$
The proof is complete.\hfill$\blacksquare$
\par
{\bf Lemma 3.4}\quad{\it  Let $R$ be a positive random variable such
that $l_F$ satisfies $SC<1$. Then
\begin{equation}
\mathcal{M}(\omega)=\{\xi+h(\omega, \xi)|\quad \xi\in E_u\}
\end{equation}
is a local invariant manifold for Equation (2.6). }
\par
{\bf Theorem 3.1 (Existence of local random invariant manifold)}
\par {\it Let $R$ be a positive random variable such that $l_F$
satisfies $SC<1$ as in the inequality (\ref{SC}). Then
\begin{equation}
\widetilde{\mathcal{M}}(\omega)=T^{-1}\mathcal{M}(\omega)
=\{\xi+e^{z(\omega)}h(\omega, e^{-z(\omega)}\xi)|\quad \xi\in E_u\}
\end{equation}
is a local invariant manifold for Equation (1.1). Namely, the graph of $
e^{z(\omega)}h(\omega, e^{-z(\omega)}\xi)$ is   the
local random invariant manifold $\widetilde{\mathcal{M}}(\omega)$ for
Equation (1.1). }
\par
Lemma 3.2, Lemma 3.4 and Theorem 3.1 can be proved as in Duan, Lu and Schmalfuss \cite{DLS04}.
\par

\renewcommand{\theequation}{\thesection.\arabic{equation}}
\setcounter{equation}{0}

\section{Local geometric shape of invariant manifolds}

\quad\quad In this section, we approximate the random invariant
manifold $M (\omega)$ step by step  to derive the local geometric
shape of the invariant manifold, as inspired by Blomker and Wang
\cite{BW}.
\par
Define
\begin{equation}
\hbar_1(t)=\int_{-\infty}^t
e^{-L_s(t-r)+\int_r^tz(\tau)d\tau}e^{-z(r)}F_s^{(R)}(e^{z(r)}v_u(r))
dr.
\end{equation}
\par
{\bf Lemma 4.1}\quad{\it There exists a positive constant $C$ such
that
\begin{equation}
\|v_s^*(t)-\hbar_1(t)\|_{C_\beta^-}\leq C|\xi|_\alpha, \quad \hbox{for arbitray }\quad t\leq 0.
\end{equation} }
\par
{\bf Proof.}\quad Note that $v_s^*=P_sv^*$. Then it follows from
Lemma 3.3 that
$$
\begin{array}{ll}
&\|v_s^*(t)-\hbar_1(t)\|_{C_\beta^-}\\
=&\|\int_{-\infty}^t
e^{-L_s(t-r)+\int_r^tz(\tau)d\tau}[e^{-z(r)}F_s^{(R)}(e^{z(r)}v^*(r))-e^{-z(r)}F_s^{(R)}(e^{z(r)}v_u^*(r))]
dr\|_{C_\beta^-}\\
\leq& \int_{-\infty}^t e^{\beta t-\int_0^tz(\tau)d\tau}
\frac{M}{(t-r)^\alpha}e^{-\lambda_s(t-r)}e^{\int_r^tz(\tau)d\tau}e^{-\beta
r
+\int_0^rz(\tau)d\tau}l_F\|v^*(r)-v_u^*(r)\|_{C_\beta^-} dr\\
=&
Ml_F\|v_s^*(r)\|_{C_\beta^-}\int_{-\infty}^t\frac{1}{(t-r)^\alpha}e^{(\beta-\lambda_s)(t-r)}dr\\
=&
Ml_F\|v_s^*(r)\|_{C_\beta^-}\frac{\Gamma(1-\alpha)}{(\lambda_s-\beta)^{1-\alpha}}\\
\leq& C|\xi|_\alpha.
\end{array}
$$
The proof is complete.\hfill$\blacksquare$
\par
{\bf Lemma 4.2}\quad{\it There exists a positive constant $C$ such
that
\begin{equation}
\|v_u^*(t)-e^{-L_ut+\int_0^tz(\tau) d\tau}\xi\|_{C_\beta^-}\leq
C|\xi|_\alpha,\quad \hbox{for arbitrary}\quad t\leq 0.
\end{equation} }
\par
{\bf Proof.}\quad Firstly, we note that
\begin{equation}
v_u^*=P_uv^*=e^{-L_ut+\int_0^tz(\tau)
d\tau}\xi+\int_0^te^{-L_u(t-r)+\int_r^tz(\tau)d\tau}e^{-z(r)}F_u^{(R)}(e^{z(r)}v(r))dr.
\end{equation}
Then for $t\leq 0$, from Lemma 3.3, we have
$$
\begin{array}{ll}
\|v_u^*(t)-e^{-L_ut+\int_0^tz(\tau) d\tau}\xi\|_{C_\beta^-}&=
\|\int_0^te^{-L_u(t-r)+\int_r^tz(\tau)d\tau}e^{-z(r)}F_u^{(R)}(e^{z(r)}v(r))dr\|_{C_\beta^-}\\
&\leq \int_0^te^{\beta
t-\int_0^tz(\tau)d\tau}e^{-\lambda_u(t-r)}e^{\int_r^tz(\tau)d\tau}
e^{-\beta
r+\int_0^rz(\tau)d\tau}l_F\|v(r)\|_{C_\beta^-}dr\\
&\leq Ml_FC|\xi|_\alpha\int_0^te^{\beta t}e^{-\lambda_u(t-r)}
e^{-\beta
r}dr\\
&\leq Ml_FC\cdot\frac{1}{\beta-\lambda_u}\cdot|\xi|_\alpha\\
&\leq C|\xi|_\alpha.
\end{array}
$$
This completes the proof.\hfill$\blacksquare$
\par

Define
\begin{equation}
\hbar_2=\int_{-\infty}^0
e^{L_sr+\int_r^0z(\tau)d\tau}e^{-z(r)}F_s^{(R)}(e^{z(r)}e^{-L_ur+\int_0^rz(\tau)
d\tau}\xi) dr.
\end{equation}
\par
{\bf Lemma 4.3}\quad{\it There exists a positive constant $C$ such
that
\begin{equation}
\|\hbar_1(0)-\hbar_2\|\leq C|\xi|_\alpha.
\end{equation} }
\par
{\bf Proof.}\quad It follows from (4.1), (4.5) and Lemma 4.2 that
$$
\begin{array}{ll}
\|\hbar_1(0)-\hbar_2\|&=\|\int_{-\infty}^0
e^{L_sr+\int_r^0z(\tau)d\tau}[e^{-z(r)}F_s^{(R)}(e^{z(r)}v_u^*(r))-e^{-z(r)}F_s^{(R)}(e^{z(r)}e^{-L_ur+\int_0^rz(\tau)
d\tau}\xi)] dr\|\\
&\leq \int_{-\infty}^0\frac{M}{(-r)^\alpha}e^{-\lambda_s(-r)}
e^{\int_r^0z(\tau)d\tau}e^{-\beta
r+\int_0^rz(\tau)d\tau}l_F\|v_u^*(r)-e^{-L_ur+\int_0^rz(\tau)
d\tau}\xi\|_{C_\beta^-} dr\\
&\leq Ml_FC|\xi|_\alpha
\int_{-\infty}^0\frac{1}{(-r)^\alpha}e^{(\lambda_s-\beta)r} dr\\
&\leq Ml_FC \cdot
\frac{\Gamma(1-\alpha)}{(\lambda_s-\beta)^{1-\alpha}}|\xi|_\alpha\\
&\leq C|\xi|_\alpha.
\end{array}
$$
The proof is thus complete.\hfill$\blacksquare$
\par
{\bf Lemma 4.4}$^{\cite{BW}}$\quad{\it There is a random variable
$K_1(\omega)$ such that $K_1(\omega)-1$ has a standard exponential
distribution and
$$
\int_0^tz(\tau)d\tau+z(t)=z(0)+\sigma \omega(t)\leq \sigma
(K_1(\omega)+|t|),\quad \hbox{for arbitrary}\quad t\leq 0,
$$
where $z$ satisfies Equation (2.4). Also,
$$
|\omega(t)|\leq \max\{\omega(t), -\omega(t)\}\leq
K^{\pm}(\omega)+|t|, \quad \hbox{for arbitrary}\quad t\leq 0,
$$
where $K^{\pm}(\omega)=K_1(\omega)+K_1(-\omega)$ and $K_1(-\omega)$ has
the same law as $K_1(\omega)$. Furthermore, for $|z(0)|$ a similar estimate is true. }
\par
Define
$$
K_2(\omega)=\sup\limits_{\tau\leq
0}|\frac{1-e^{-\lambda_u\tau+\sigma \omega(\tau)}}{\gamma e^{\delta
|\tau|}}|,
$$
where $\gamma$ and $\delta$ are the positive constants.
\par
{\bf Lemma 4.5}\quad {\it Choose two positive real parameters
$\gamma$ and $\delta$ satisfying $\gamma\geq \max\{-\lambda_u,
\sigma\}$ and  $\delta>-\lambda_u+\sigma$. Then there is a constant
$C$ such that
$$
K_2(\omega)\leq Ce^{\sigma K^{\pm}(\omega)}(1+K^{\pm}(\omega)).
$$
}
\par
{\bf Proof.}\quad Using $|1-e^x|\leq |x|e^{|x|}$ and Lemma 4.4, we
get
$$
\begin{array}{ll}
K_2(\omega)&=\sup\limits_{\tau\leq
0}|\frac{1-e^{-\lambda_u\tau+\sigma \omega(\tau)}}{\gamma e^{\delta
|\tau|}}|\\
& \leq \sup\limits_{\tau\leq
0}\frac{|\lambda_u||\tau|+|\sigma||\omega(\tau)|}{\gamma}e^{|\lambda_u||\tau|+|\sigma||\omega(\tau)|}e^{-\delta
|\tau|}\\
& \leq e^{\sigma K^{\pm}(\omega)}\sup\limits_{\tau\leq
0}(|\tau|+|\omega(\tau)|)e^{(-\lambda_u+\sigma-\delta)|\tau|}\\
& \leq Ce^{\sigma K^{\pm}(\omega)}\sup\limits_{\tau\leq
0}(|\tau|e^{(-\lambda_u+\sigma-\delta)|\tau|}+K^{\pm}(\omega)e^{(-\lambda_u+\sigma-\delta)|\tau|})\\
& \leq Ce^{\sigma K^{\pm}(\omega)}(1+K^{\pm}(\omega)).
\end{array}
$$
The proof is complete.\hfill$\blacksquare$
\par
{\bf Lemma 4.6}\quad{\it Let $e^{z(0)}|\xi|_\alpha\leq R$. Then
there exists a positive constant $C$ such that
\begin{equation}
|\chi_R(e^{z(r)}e^{-L_ur+\int_0^rz(\tau) d\tau}\xi)-1|\leq
\frac{C}{R}e^{z(0)}K_2(\omega)\gamma e^{-\delta r}|\xi|_\alpha,\quad \hbox{for arbitrary}\quad r\leq 0.
\end{equation} }
\par
{\bf Proof.}\quad Note that $e^{z(0)}|\xi|_\alpha\leq R$. Then
$\chi_R(e^{z(0)}\xi)=1$. Therefore, for $r\leq 0$, it follows from
Lemma 4.4 that
\begin{equation}
\begin{array}{ll}
|\chi_R(e^{z(r)}e^{-L_ur+\int_0^rz(\tau) d\tau}\xi)-1|&\leq
|\chi_R(e^{z(r)}e^{-\lambda_ur+\int_0^rz(\tau)
d\tau}\xi)-\chi_R(e^{z(0)}\xi)|\\
&\leq \frac{C}{R}|e^{-\lambda_ur+z(r)+\int_0^rz(\tau)
d\tau}\xi-e^{z(0)}\xi|_\alpha\\
&\leq \frac{C}{R}|\xi|_\alpha\cdot
|e^{-\lambda_ur+z(0)+\sigma \omega(r)}-e^{z(0)}|\\
&\leq \frac{C}{R}e^{z(0)}|\xi|_\alpha \cdot
|1-e^{-\lambda_ur+\sigma \omega(r)}|\\
&\leq \frac{C}{R}e^{z(0)}K_2(\omega)\gamma e^{-\delta
r}|\xi|_\alpha.
\end{array}
\end{equation}
The proof is complete.\hfill$\blacksquare$
\par

We further define
\begin{equation}
\hbar_3=\int_{-\infty}^0
e^{L_sr+\int_r^0z(\tau)d\tau}e^{-z(r)}F_s(e^{z(r)}e^{-L_ur+\int_0^rz(\tau)
d\tau}\xi) dr.
\end{equation}
\par
{\bf Lemma 4.7}\quad{\it Let
$0<\sigma<\frac{\lambda_s-(p-1)\lambda_u}{p}$ and
$e^{z(0)}|\xi|_\alpha\leq R$. Then there exists a positive constant
$C$ such that
\begin{equation}
\|\hbar_2-\hbar_3\|\leq Ce^{z(0)}K_2(\omega)e^{(p-1)\sigma
K_1(\omega)}|\xi|_\alpha^{2}.
\end{equation} }
\par
{\bf Proof.}\quad Firstly, from the condition $\sigma<\frac{\lambda_s-(p-1)\lambda_u}{p}$,
we know that
$$
\lambda_s-p\lambda_u-(p-1)\sigma> -\lambda_u+\sigma,
$$
which implies that there must exist a constant $\delta$ satisfying
\begin{equation}
\lambda_s-p\lambda_u-(p-1)\sigma>\delta> -\lambda_u+\sigma.
\end{equation}
\par
Also, we note that $F(u)=u^p$, $F^{(R)}(u)=\chi_R(u)F(u)$,
$e^{z(0)}|\xi|_\alpha\leq R$ and (2.7). Therefore, from Lemma 4.7
and Lemma 4.4, we get
$$
\begin{array}{ll}
&\|\hbar_2-\hbar_3\|\\
\leq &\|\int_{-\infty}^0
e^{L_sr+\int_r^0z(\tau)d\tau}e^{-z(r)}[F_s^{(R)}(e^{z(r)}e^{-L_ur+\int_0^rz(\tau)
d\tau}\xi)-F_s(e^{z(r)}e^{-L_ur+\int_0^rz(\tau) d\tau}\xi)]
dr\|\\
\leq &\|\int_{-\infty}^0
e^{L_sr+\int_r^0z(\tau)d\tau}e^{(p-1)z(r)}e^{-pL_ur+p\int_0^rz(\tau)
d\tau}\xi_s^p[\chi_R(e^{z(r)}e^{-L_ur+\int_0^rz(\tau) d\tau}\xi)-1]
dr\|\\
\leq & \frac{C}{R}e^{-(p-1)z(0)}K_2(\omega)\gamma |\xi|_\alpha
\int_{-\infty}^0e^{-\delta r}
e^{L_sr+\int_r^0z(\tau)d\tau}e^{(p-1)z(r)}e^{-pL_ur+p\int_0^rz(\tau)
d\tau}\|(e^{z(0)}\xi_s)^p\|
dr\\
\leq & \frac{C}{R}e^{-(p-2)z(0)}K_2(\omega)\gamma
|\xi|_\alpha^{2}M^2l_F\int_{-\infty}^0\frac{1}{(-r)^\alpha}
e^{(p-1)[z(r)+\int_0^rz(\tau)d\tau]}e^{(\lambda_s-p\lambda_u-\delta)r}dr\\
\leq & \frac{C}{R}e^{-(p-2)z(0)}K_2(\omega)\gamma
|\xi|_\alpha^{2}M^2l_Fe^{(p-1)\sigma
K_1(\omega)}\int_{-\infty}^0\frac{1}{(-r)^\alpha}
e^{(\lambda_s-p\lambda_u-\delta-(p-1)\sigma)r}dr,
\end{array}
$$
which, from (4.11), immediately implies
that
$$
\begin{array}{ll}
\|\hbar_2-\hbar_3\| &\leq \frac{C}{R}e^{-(p-2)z(0)}K_2(\omega)\gamma
|\xi|_\alpha^{p+1}M^2l_Fe^{(p-1)\sigma
K_1(\omega)}\frac{\Gamma(1-\alpha)}{(\lambda_s-p\lambda_u-\delta-(p-1)\sigma)^{1-\alpha}}\\
&\leq CK_2(\omega)e^{(p-1)\sigma
K_1(\omega)-(p-2)z(0)}|\xi|_\alpha^{2}.
\end{array}
$$
The proof is complete.\hfill$\blacksquare$
\par
{\bf Lemma 4.8}\quad {\it Let
\begin{equation}
K_3(\omega)=\sup\limits_{r\leq
0}|\frac{1-e^{(p-1)\sigma\omega(r)}}{\gamma_1e^{(p-1)\delta_1|r|}}|.
\end{equation}
If $\gamma_1>\sigma$ and $\delta_1>\sigma$, then
$$
K_3(\omega)\leq
Ce^{(p-1)\sigma K^{\pm}(\omega)}(1+K^{\pm}(\omega)).
$$
Furthermore,
\begin{equation}
|1-e^{(p-1)\sigma\omega(r)}|\leq
K_3(\omega)\gamma_1e^{-(p-1)\delta_1r},\quad \hbox{for arbitrary}\quad r\leq 0.
\end{equation}
}
\par
{\bf Proof.}\quad Using $|1-e^x|\leq |x|e^{|x|}$ and Lemma 4.4, we
get
$$
\begin{array}{ll}
K_3(\omega)&=\sup\limits_{r\leq
0}|\frac{1-e^{(p-1)\sigma\omega(r)}}{\gamma_1e^{(p-1)\delta_1|r|}}|\\
& \leq \sup\limits_{r\leq
0}\frac{(p-1)\sigma|\omega(r)|}{\gamma_1}e^{(p-1)\sigma|\omega(r)|}e^{-(p-1)\delta_1|r|}\\
& \leq (p-1)e^{(p-1)\sigma K^{\pm}(\omega)}\sup\limits_{r\leq
0}(|r|+K^\pm(\omega))e^{(p-1)(\sigma-\delta_1)|r|}\\
& \leq (p-1)e^{(p-1)\sigma K^{\pm}(\omega)}\sup\limits_{r\leq
0}(|r|e^{(p-1)(\sigma-\delta_1)|r|}+K^\pm(\omega)e^{(p-1)(\sigma-\delta_1)|r|})\\
& \leq Ce^{\sigma K^{\pm}(\omega)}(1+K^{\pm}(\omega)).
\end{array}
$$
The proof is complete.\hfill$\blacksquare$

\par
{\bf Lemma 4.9}\quad{\it Let $0<\sigma<-\lambda_u$ and
$e^{z(0)}|\xi|_\alpha\leq R$. Then there exists a positive constant
$C$ such that
\begin{equation}
\|\hbar_3-e^{(p-1)z(0)}(L_s-pL_u)^{-1}\xi_s^p\|\leq  Ce^{(p-1)z(0)} K_3(\omega) |\xi|_\alpha.
\end{equation} }
\par
{\bf Proof.}\quad It follows from (4.9) and Lemma 4.4 that
$$
\hbar_3=e^{(p-1)z(0)}\xi_s^p\int_{-\infty}^0 e^{L_s r} e^{-pL_ur}
e^{(p-1)\sigma\omega(r)} dr.
$$
\par
From the condition $\sigma<-\lambda_u$, there must exist a parameter
$\delta_1$ satisfying $\sigma<\delta_1<-\lambda_u$, which implies
that
\begin{equation}
\lambda_s-p\lambda_u-(p-1)\delta_1>0.
\end{equation}
\par
Also note that
$$
e^{(p-1)z(0)}\xi_s^p\int_{-\infty}^0 e^{L_s r} e^{-pL_u r}
dr=e^{(p-1)z(0)}(L_s-pL_u)^{-1}\xi_s^p .
$$
\par
Therefore, it follows from Lemma 4.8 and (4.15) that
$$
\begin{array}{ll}
&\|\hbar_3-e^{(p-1)z(0)}(L_s-pL_u)^{-1}\xi_s^p\|\\
=& \|e^{-z(0)}(e^{z(0)}\xi_s)^p\int_{-\infty}^0 e^{L_s r} e^{-pL_ur}
[e^{(p-1)\sigma\omega(r)}-1] dr\|\\
\leq & |\xi|_\alpha K_3(\omega)\gamma_1 M^2l_F
\int_{-\infty}^0
\frac{1}{(-r)^\alpha}e^{(\lambda_s-p\lambda_u-(p-1)\delta_1)r}dr\\
\leq &|\xi|_\alpha K_3(\omega)\gamma_1 M^2l_F
\frac{\Gamma(1-\alpha)}{(\lambda_s-p\lambda_u-(p-1)\delta_1)^{1-\alpha}}\\
\leq & CK_3(\omega) |\xi|_\alpha .\\
\end{array}
$$
The proof is complete.\hfill$\blacksquare$
\par
{\bf Lemma 4.10}\quad{\it Let
$0<\sigma<\min\{\frac{\lambda_s-(p-1)\lambda_u}{p}, -\lambda_u\}$
and $e^{z(0)}|\xi|_\alpha\leq R$. Then there exists a positive
constant $C$ such that
\begin{equation}
\|h(\omega, \xi)-e^{(p-1)z(0)}(L_s-pL_u)^{-1}\xi_s^p\|\leq C[1
+K_3(\omega)+e^{(p-1)\sigma K_1(\omega)-(p-2)z(0)}K_2(\omega)|\xi|_\alpha]\cdot|\xi|_\alpha.
\end{equation} }
\par
This lemma is directly   from Lemma 4.1, Lemma 4.3, Lemma
4.7 and Lemma 4.9.
\par
Finally we have the following main result.
\par
{\bf Theorem 4.1 (Local geometric shape of random invariant manifold)}
\par {\it Let $0<\sigma<\min\{\frac{\lambda_s-(p-1)\lambda_u}{p},
-\lambda_u\}$ and $|\xi|_\alpha\leq R$. Then there exists a positive
constant $C$ such that
\begin{equation}
\|e^{z(\omega)}h(\omega,
e^{-z(\omega)}\xi)-(L_s-pL_u)^{-1}\xi_s^p\|\leq
C(|\xi|_\alpha+|\xi|_\alpha^2)
\end{equation}
holds with probability larger than $1-Ce^{-\frac{1}{\sigma}}$.
Therefore in a neighborhood of zero for Equation (1.1), the graph
$(\xi, e^{z(\omega)}h(\omega, e^{-z(\omega)}\xi))$ of the invariant
manifold $\widetilde{\mathcal{M}}(\omega)$ is approximately given by
$(\xi, (L_s-pL_u)^{-1}\xi_s^p)$ with probability larger than
$1-Ce^{-\frac{1}{\sigma}}$. }

\par
{\bf Proof.}\quad Define $\Omega_K=\{\omega\in \Omega|\quad
K^\pm(\omega)>\frac{1}{\sigma}\}$. By Lemma 4.4 this set has
probability less than $Ce^{-\frac{1}{\sigma}}$. Therefore, on the
complement $\Omega_K^C$, there must exist a positive constant $C$
such that
$$
K_1(\omega)\leq C, K_2(\omega)\leq C, K_3(\omega)\leq C.
$$
Therefore, it follows from Lemma 4.10 that Theorem 4.1 holds. The
proof is complete.\hfill$\blacksquare$
\par
{\bf Remark 4.1}\quad {\it Here we present an example to explain
Theorem 4.1. Consider Equation (1.1) with the line operator
$L=-\partial_{xx}-3\cdot id$ on $[0, \pi]$ with the homogeneous
Dirichlet boundary condition. Then the eigenvalues of $L$ are
$$
\lambda_1=-2,\;\lambda_2=1,\;\lambda_3=6,\;\cdots,
\lambda_k=k^2-3,\;\cdots
$$
with the corresponding eigenfunctions
$$
e_1=\sin x,\; e_2=\sin 2x, \;e_3=\sin 3x,\;\cdots, e_k=\sin kx,\;
\cdots
$$
Therefore, Theorem 4.1 affords a local unstable invariant manifold
for Equation (1.1). In this case, $L_u=-2\cdot id$ and
$E_u=span\{e_1\}$. Then we can write $\xi=r\cdot e_1 $ with $r\in
\mathbb{R}$. We denote $e_1^{\bot}:=(L_s+2p\cdot id)^{-1}P_se_1^p$.
Then
$$
(L_s-pL_u)^{-1}\xi_s^p=r^p\cdot(L_s+2p\cdot
id)^{-1}P_se_1^p=r^p\cdot e_1^{\bot}.
$$
Therefore in the state space spanned by the coordinate variables
$e_1$ and $e_1^{\bot}$, the geometric shape of
$\widetilde{\mathcal{M}}(\omega)$ is given by $(r, r^p)$.}
\par

\renewcommand{\theequation}{\thesection.\arabic{equation}}
\setcounter{equation}{0}

\section{Results for the corresponding deterministic system}

\quad\quad In this section, we briefly comment on invariant
manifolds for the corresponding deterministic system of Equation
(1.1)(i.e. Equation (1.1) with $\sigma=0$). We consider the local
(unstable) invariant manifold and its local geometric shape, for
this deterministic system. Since we use the same method as in
Section 3 and in Section 4 above, we omit the proofs of the results
but only highlight some differences.
\par
Consider the deterministic system
\begin{equation}
\frac{d u}{d t}+Lu-u^p=0.
\end{equation}
with the initial data $u(x, 0)=u_0=x\in E^\alpha$.
\par
Define a Banach space for each $\beta\in (\lambda_u, \lambda_s)$ as
follows
$$
\mathfrak{C}_{\beta}^-=\{f(\cdot)\in C((-\infty, 0]; E^\alpha)|\quad
\sup\limits_{t\leq 0}e^{\beta t} |f|_\alpha<\infty \}
$$
with the norm
$$
\|f\|_{\mathfrak{C}_{\beta}^-}=\sup\limits_{t\leq 0}e^{\beta t}
|f|_\alpha.
$$
\par
{\bf Lemma 5.1}\quad{\it Assume that $u(\cdot)$ is in
$\mathfrak{C}_{\beta}^-$. Then $u(t)$ is the local solution of
Equation (5.1) with the initial datum $u(0)=x$ if and only if $u(t)$
satisfies
\begin{equation}
u(t)=e^{-L_ut}\xi+\int_0^te^{-L_u(t-r)}F_u^{(R)}(u(r))dr
+\int_{-\infty}^t e^{-L_s(t-r)}F_s^{(R)}(u(r)) dr,
\end{equation}
where $\xi=P_ux\in E_u$. }
\par
Let $\mathfrak{J}(u, \xi)$ denote the right hand side of (5.2), and
$h(\xi)=\int_{-\infty}^0 e^{-L_s(t-r)}F_s^{(R)}(u(r)) dr$. Then we
have following results.
\par
{\bf Lemma 5.2}\quad{\it Let $R$ be a positive real number such that
$l_F$ satisfies $SC<1$. Then $\mathfrak{J}(u, \xi)$ has a unique
fixed point $u^*=u^*(t; \xi)=\mathfrak{J}(u^*)\in
\mathfrak{C}_\beta^-$. Furthermore, there exist a positive constant
$C$ such that
$$
\|u^*(t;\xi_1)-u^*(t;\xi_2)\|_{\mathfrak{C}_\beta^-}\leq
C|\xi_1-\xi_2|_\alpha,
$$
$$ \|u^*(t;\xi)\|_{\mathfrak{C}_\beta^-}\leq C|\xi|_\alpha,\quad
\|u_s^*(t;\xi)\|_{\mathfrak{C}_\beta^-}\leq C|\xi|_\alpha,\quad
and\quad \|u_u^*(t;\xi)\|_{\mathfrak{C}_\beta^-}\leq C|\xi|_\alpha,
$$
where $u_s^*=P_su^*$ and $u_u^*=P_uu^*$. Moreover,
$$
\|h(\xi_1)-h(\xi_2)\|_{\mathfrak{C}_\beta^-}\leq
C|\xi_1-\xi_2|_\alpha.
$$. }
\par
{\bf Theorem 5.1 (Existence of local invariant manifold for the
corresponding deterministic system)}\par {\it Let $R$ be a positive
real number such that $l_F$ satisfies $SC<1$ as in (\ref{SC}). Then
\begin{equation}
\widetilde{\mathcal{M}}=\{\xi+h(\xi)|\quad \xi\in E_u\}
\end{equation}
is a local invariant manifold for the deterministic Equation (5.1).  Namely, the graph of  $
h(\xi)$ is  the local  deterministic   invariant manifold
$\widetilde{\mathcal{M}}$ for Equation (5.1).}
\par
In the following, we approximate the local invariant manifold
$\widetilde{\mathcal{M}}$  for Equation (5.1). Define
\begin{equation}
\hbar_1(t)=\int_{-\infty}^t e^{-L_s(t-r)}F_s^{(R)}(u_u(r)) dr,
\end{equation}
\begin{equation}
\hbar_2=\int_{-\infty}^0 e^{L_sr}F_s^{(R)}(e^{-L_ur}\xi) dr,
\end{equation}
and
\begin{equation}
\hbar_3=\int_{-\infty}^0 e^{L_sr}F_s(e^{-L_ur}\xi) dr.
\end{equation}
\par
{\bf Lemma 5.3}\quad{\it Let $|\xi|_\alpha\leq R$. There exists a
positive constant $C$ such that
$$
\|u_s^*(t)-\hbar_1(t)\|_{\mathfrak{C}_\beta^-}\leq C|\xi|_\alpha,
\quad \hbox{for arbitray }\quad t\leq 0,
$$
$$
\|u_u^*(t)-e^{-L_ut+\int_0^tz(\tau)
d\tau}\xi\|_{\mathfrak{C}_\beta^-}\leq C|\xi|_\alpha,\quad \hbox{for
arbitrary}\quad t\leq 0,
$$
$$
|\chi_R(e^{-L_ut}\xi)-1|\leq
\frac{C}{R}(1-e^{-\lambda_ut})|\xi|_\alpha,\quad \hbox{for
arbitrary}\quad t\leq 0,
$$
and
$$
\|\hbar_1(0)-\hbar_2\|\leq C|\xi|_\alpha,\quad
\|\hbar_2-\hbar_3\|\leq C|\xi|_\alpha^{2}.
$$
}
\par
Furthermore, noting that $\hbar_3=(L_s-pL_u)^{-1}\xi_s^p$ and Lemma
5.3, we have the following result about local geometric shape of
invariant manifold for Equation (5.1).
\par
{\bf Theorem 5.2 (Local geometric shape of invariant manifold for the
corresponding deterministic system)}
\par{\it Let $|\xi|_\alpha\leq R$. Then there exists a positive
constant $C$ such that
\begin{equation}
\|h(\xi)-(L_s-pL_u)^{-1}\xi_s^p\|\leq
C(|\xi|_\alpha+|\xi|_\alpha^2).
\end{equation}
Therefore in a neighborhood of zero for Equation (5.1), the graph
$(\xi, h(\xi))$ of the invariant manifold $\widetilde{\mathcal{M}}$
is approximately given by $(\xi, (L_s-pL_u)^{-1}\xi_s^p)$.}
\par

\renewcommand{\theequation}{\thesection.\arabic{equation}}
\setcounter{equation}{0}

\section{Conclusions}

\quad\quad For a class of stochastic partial
differential equations,  after establishing the existence of  the local unstable random invariant manifold
(see Theorem 3.1), we derive an approximation for the local geometric shape of this random invariant
manifold (see Theorem 4.1). The local geometric shape approximation holds with
significant probability. Furthermore, with the  noise  intensity
  $\sigma$ decreasing, this significant  probability   is increasing.
In fact, as noise intensity $\sigma\searrow 0$, the probability
$1-Ce^{-\frac{1}{\sigma}} \nearrow 1$. On the other hand, when
$\sigma=0$, Equation (1.1) is a deterministic system, the local
geometric shape approximation of the corresponding deterministic
invariant manifold is the same but holds surely (see Theorem 5.2).



\begin{thebibliography}{99}



\bibitem{Arnold}
L. Arnold, Random Dynamical Systems, Springer-Verlag, 1998.

\bibitem{B}
D. Blomker, Amplitude Equations for Stochastic Partial Differential Equations, vol. 3 of Interdisciplinary
Mathematical Sciences. World Scientific Publishing, Singapore, 2007.

\bibitem{BW}
D. Blomker and W. Wang, Qualitative properties of local random
invariant manifolds for SPDE with quadratic nonlinearity, \emph{ J.
Dyn. Diff. Equat.}, 2009, DOI 10.1007/s10884-009-9145-6.

\bibitem{CDKS}  T. Caraballo, J. Duan, K. Lu and B. Schmalfuss,  Invariant
manifolds for random and stochastic partial differential equations.
\emph{Advanced Nonlinear Studies}, \textbf{10} (2009): 23-52.

\bibitem{CCL05}
T. Caraballo, I. Chueshov and J. Langa, Existence of invariant
manifolds for coupled parabolic and hyperbolic stochastic partial
differential equations, \emph{Nonlinearity}, \textbf{18} (2) (2005):
747-767.

\bibitem{DZ}
G. Da Prato and J. Zabczyk, Stochastic Equations in Infinite
Dimensions. Cambridge University Press, Cambridge, 1992.

\bibitem{DLS03}
J. Duan, K. Lu and B. Schmalfuss, Invariant manifolds for stochastic
partial differential equations, \emph{ Ann. Probab.}, \textbf{31}
(4) (2003): 2109-2135.

\bibitem{DLS04}
J. Duan, K. Lu and B. Schmalfuss, Smooth stable and unstable
manifolds for stochastic evolutionary equations, \emph{J. Dynam.
Differential Equations}, \textbf{16} (4) (2004): 949-972.

\bibitem{GS}
J. Garcia-Ojalvo and J. M. Sancho, Noise in Spatially Extended
Systems. Springer, Berlin, 1999.

\bibitem{Henry}
D. Henry, Geometric Theory of Semilinear Parabolic Equations, Lecture Notes in
Mathematics, Vol. 840, Springer-Verlag, 1981.

\bibitem{LS07}
K. Lu and B. Schmalfuss, Invariant manifolds for stochastic wave
equations, \emph{J. Differential Equations}, \textbf{236} (2007):
460-492.

\bibitem{MZZ08}
S. A. Mohammed, T. Zhang and H. Zhao, The stable manifold theorem
for semilinear stochastic evolution equations and stochastic partial
differential equations, \emph{ Memoirs of the American Mathematical
Society}, \textbf{196} (2008): 1-105.


 \bibitem{Rockner} C. Prevot and M. Rockner,
\emph{A Concise Course on Stochastic Partial Differential
Equations}, Lecture Notes in Mathematics, Vol. 1905. Springer, New
York, 2007.

\bibitem{R}
B. L. Rozovskii, Stochastic Evolution Equations. Kluwer, Boston,
1990.

\bibitem{SDL}
X. Sun, J. Duan and X. Li, An impact of noise on invariant manifolds
in nonlinear dynamical systems, \emph{J. Math. Phys.}, \textbf{51},
042702 (2010).

\bibitem{T}
R. Temam, Infinite Dimensional Dynamical Systems in
Mechanics and Physics. Springer-Verlag, New York, 1997.

\bibitem{WD}
E. Waymire and J. Duan, Probability and Partial Differential
Equations in Modern Applied Mathematics, IMA vol. 140. Springer, New
York, 2005.








\end{thebibliography}
\end{document}